\documentclass[conference, 9pt]{IEEEtran}
\IEEEoverridecommandlockouts
\usepackage{cite}
\usepackage{amsmath,amssymb,amsfonts}
\usepackage{graphicx}
\usepackage{textcomp}
\usepackage{xcolor}
\usepackage{booktabs}
\usepackage{tikz}
\usepackage{color}
\usepackage{graphicx}
\usepackage{subfigure} 
\usepackage{float}
\usepackage{dsfont}
\usepackage{algorithm}
\usepackage{algpseudocode}
\usepackage{amsfonts}

\usepackage{amssymb}%
\usepackage{multirow} 
\usepackage{makecell}
\usepackage{diagbox}
\usepackage{threeparttable}
\usepackage{makecell}
\usepackage{balance}

\usepackage{changes}
\usepackage{CJKutf8}

\def\BibTeX{{\rm B\kern-.05em{\sc i\kern-.025em b}\kern-.08em
    T\kern-.1667em\lower.7ex\hbox{E}\kern-.125emX}}
\begin{document}
\begin{CJK*}{UTF8}{gbsn}

\title{Per-RMAP: Feasibility-Seeking and Superiorization Methods for Floorplanning with I/O Assignment
}

\author{\IEEEauthorblockN{Shan Yu\textsuperscript{1}, Yair Censor\textsuperscript{2}, Ming Jiang\textsuperscript{1} and Guojie Luo\textsuperscript{3,4}}
\IEEEauthorblockA{\textit{\textsuperscript{1}Department of Information and Computational Sciences, School of Mathematical Sciences, Peking University, Beijing, China} \\
\textit{\textsuperscript{2}Department of Mathematics, University of Haifa, Haifa, Israel}\\
\textit{\textsuperscript{3}National Key Laboratory for Multimedia Information Processing, School of Computer Science, Peking University, Beijing, China}\\
\textit{\textsuperscript{4}Center for Energy-efficient Computing and Applications, Peking University, Beijing, China}\\
\textit{yu-shan@pku.edu.cn, yair@math.haifa.ac.il, \{ming-jiang, gluo\}@pku.edu.cn}\\\textit{March 31, 2023}
}
}

\maketitle

\begin{abstract}
The feasibility-seeking approach provides a systematic scheme to manage and solve complex constraints for continuous problems, and we explore it for the floorplanning problems with increasingly heterogeneous constraints. The classic legality constraints can be formulated as the union of convex sets. However, the convergence of conventional projection-based algorithms is not guaranteed when the constraints sets are non-convex, which is the case with unions of convex sets. In this work, we propose a resetting strategy to greatly eliminate the divergence issue of the projection-based algorithm for the feasibility-seeking formulation. Furthermore, the superiorization methodology (SM), which lies between feasibility-seeking and constrained optimization, is firstly applied to floorplanning. The SM uses perturbations to steer the iterates of a feasibility-seeking algorithm to a feasible solution with shorter total wirelength. The proposed algorithmic flow is extendable to tackle various constraints and variants of floorplanning problems, e.g., floorplanning with I/O assignment problems. We have evaluated the proposed algorithm on the MCNC benchmarks. We can obtain legal floorplans only two times slower than the branch-and-bound method in its current prototype using MATLAB, with only 3\% wirelength inferior to the optimal results. We evaluate the effectiveness of the algorithmic flow by considering the constraints of I/O assignment, and our algorithm achieves 8\% improvement on wirelength.
\end{abstract}

\begin{IEEEkeywords}
feasibility-seeking, superiorization method, projection algorithms,
floorplanning, I/O assignment.
\end{IEEEkeywords}

\section{Introduction}
Floorplanning is a critical stage of the VLSI physical design flow,
because it influences the quality of down-stream stages. It can be
described as the task of placing a given set of rectangular modules\footnote{ The floorplanning considers two kinds of modules, hard modules
and soft modules. Soft modules have a fixed area while allowing variable
heights and widths whereas hard modules have fixed heights and widths.
In this paper, we consider only hard modules. The approach can be
extended to tackle soft modules in future work.} into a rectangular region such that there are no overlaps among modules,
while minimizing wirelength, congestion, temperature, etc. It is a
hard problem~\cite{Adya2003}, and numerous and diverse conditions
could be taken into consideration to yield more effective integrated
circuits, such as boundary conditions, non-overlap conditions and
I/O assignment conditions. The \textbf{feasibility-seeking problem
(FSP)} is to find constraint-compatible points for a family (usually
finite) of constraints sets. At the same time, as further explained
below, FSP enables to use superiorization which adopts the
philosophy of ``satisficing'' rather than optimizing when modeling
and solving a problem that includes both constraints and an objective
function.

State-of-the-art floorplanning methods can be roughly divided into
four categories: meta-heuristic methods, exact methods (e.g., branch-and-bound (B\&B) methods), analytical methods, and learning-based methods. Heuristic methods and branch-and-bound methods first adopt a representation, such as a sequence pair \cite{Murata1996} and a $B^{*}$-tree\cite{Chang2000}, and then search for the optimal or a sub-optimal solution in the representation space. Heuristic methods like \cite{Cong2004} search heuristically and stop if certain criteria are achieved whereas exact methods search in the whole search space to find the optimal solution. B\&B methods are most important methods in the class of exact methods, which build a rooted decision tree to enumerate the search space and reduce it by pruning useless branches, see, e.g., \cite{funke2016exact}. Analytical methods model the floorplanning as an optimization problem with quadratic \cite{Zhan2006} or nonlinear \cite{Lin2014} objective functions. Generally, they employ a global floorplanning step to get the rough position of all modules, and use the legalization step to eliminate the overlaps and get the exact positions of all modules. Recently learning-based algorithms, especially reinforcement learning, became popular. The GoodFloorplan~\cite{Xu2022} combines graph convolutional network and reinforcement learning to explore the design space. Liu et al~\cite{Liu2022} use graph attention to learn an optimized mapping between circuit connectivity and physical wirelength, and produce a chip floorplan using efficient model inference.

Compared with the above methods, the FSP formulation, explored in
this paper, employes some powerful tools to tackle problems with numerous
and diverse constraints. Heuristic and AI methods need delicate design
which may adversely affect diverse constraints. Exact methods may
have a huge search space. Analytic methods will find it challenging
to keep a trade-off between constraints and objectives. Compared with
those, FSP focuses on the feasibility and simplifies the problem. Additionally, iterative projection methods, frequently used tools in FSP, are fast, easy to implement and to tackle complex constraints sets in floorplanning. The algorithmic flow makes it easy to handle extended constraints, for example, in our paper I/O assignment is taken into consideration beyond basic floorplanning, and it achieves $8\%$ improvement on the total wirelength within bearable time cost.

Furthermore, the \textbf{superiorization methodology} (SM) could be utilized
to find superior solutions. It uses perturbations to steer a feasibility-seeking
algorithm to an output that is feasible and of which a certain target function is improved. A recent tutorial~\cite{Yair2022} on the SM
contains pointers to a variety of recent works and sources on this
subject. Blake Schultze et al. \cite{yair2020} applies the total
variation (TV) superiorization to image reconstruction in proton computed
tomography to improve the result. Inspired by their approach, we propose
a wirelength-superiorized FSP algorithm for floorplanning.

In short, our contributions are: 
\begin{enumerate}
\item We formulate floorplanning as a feasibility-seeking problem, which
focus on the feasibility and simplifies the problem. 
\item We propose a generalized projection method using a resetting strategy to tackle complex constraints sets in floorplanning. This resetting strategy improves the initial convergence behavior of the projection method.
\item We apply a wirelength superiorization method to find a superior feasible solution with shorter total wirelength. 
\item Our algorithmic flow shows potential to tackle diverse constraints.
With considering I/O assignment for floorplanning, our flow achieves
$8\%$ improvement compared with B\&B method\footnote{The B\&B method for comparison does not consider I/O assignment. The search space and pruning strategies will have to be completely re-defined and re-implemented when considering extra variables or constraints in the B\&B method. By contrast, the FSP method easily adapts to a new formulation.} on
the total wirelength within bearable time. 
\end{enumerate}

The remainder of the paper is organized as follows: Section~\ref{chap:formulation}
presents the feasibility-seeking formulation in the context of the
floorplanning problem. Section~\ref{chap:algorithm} describes our
proposed method, called perturbed resettable method of alternating projection (Per-RMAP). Section~\ref{chap:experiment}
provides experimental results, followed by some concluding comments
in Section~\ref{chap:conclusion}.
\section{Floorplanning as a Feasibility Problem}
\label{chap:formulation}

In this section, the floorplanning problem is formulated as a feasibility-seeking problem of finding a point in the intersection of a finite number of sets and a projection method is introduced to find solutions.

\subsection{Design Criteria for the Floorplanning Problem}
The problem of floorplanning with I/O assignment is to
find the positions of modules in a retangular floorplanning region and the positions of I/O pins on the boundary of the floorplanning region while reducing the total wirelength of connections among them. 

A solution of the the floorplanning problem should satisfy the following
conditions. 
\begin{enumerate}
\item[I] \textbf{Boundary}: Every module is within the given floorplanning region; 
\item[II] \textbf{Non-overlap}: Every pair of modules has no overlap.
\end{enumerate}

Pins are the locations for wire connections on the boundaries of modules and on the boundary of the floorplanning region, while nets denote sets of pins that require electrical connections.

The solution for the floorplanning problem posed as a feasibility-seeking
problem is in general not unique. We seek to find a solution satisfying
the following additional condition. 
\begin{enumerate}
\item[III] \textbf{Wirelength}: The total wirelength of nets in terms of the half-perimeter wirelength (HPWL) is as short as possible. 
\end{enumerate}

Let $\mathcal{P}$ be the set of the pins. The $p$-th pin $P_{p}\in\mathcal{P}$ is the point located at $(x_{p}^{pin},y_{p}^{pin})$. Let $\mathcal{E}$ denote the nets of wire connections among connected pins from $\mathcal{P}$. Then the HPWL
of a floorplanning is 

\begin{equation}
  \begin{split}
  \mbox{\textrm{HPWL}} &(x,y):=\\
  & \sum_{e\in\mathcal{E}}(\max_{p,q\in e}\vert x_{p}^{pin}-x_{q}^{pin}\vert + \max_{p,q\in e}\vert y_{p}^{pin}-y_{q}^{pin}\vert).
  \end{split}
\end{equation}

The above conditions I, II and III are our design criteria for the floorplanning problem in this work. Of course, there are other design criteria such as power distribution and module density, etc., which should be considered in practice. We believe that those extra criteria can also be formulated as feasibility constraints but they are left for future study.

\subsection{Feasibility-Seeking Formulation for Floorplanning Problem}
To reach our feasibility-seeking formulation for the floorplanning
problem, we need more notations. The floorplanning region is a 2D rectangle with the left bottom corner at the origin $(0,0)$ and the upper right
corner at $(W,H)$ where $W\:\mathrm{and\:}H$ are some given positive
real numbers. Let $\mathcal{M}$ be the set of modules to be placed
into the given floorplanning region. For the $i$-th module $m_{i}\in\mathcal{M}$,
its width and height are $w_{i}$ and $h_{i}$, respectively, its
``coordinate'' is the location of its bottom left corner, denoted
by $(x_{i},y_{i})$ which is to be determined. Assume pin $P_p\in\mathcal{P}$ belongs to module $m_i$, its location $(x_p^{pin},y_p^{pin})=(x_i+x_p^{\mathit{offset}},y_i+y_p^{\mathit{offset}})$ is determined by $(x_i,y_i)$ with constant pin offset $x_p^{\mathit{offset}}$ and $y_p^{\mathit{offset}}$ for hard modules. Let $\mathcal{P}_{io}\subset\mathcal{P}$
be the set of the I/O pins that are on the boundary of the floorplanning region. The coordinates of the I/O pins
$P_{p}\in\mathcal{P}_{io}$ at $(x_{p}^{pin},y_{p}^{pin})$ are to
be determined when considering I/O assignment. Let the total number
of modules in $\mathcal{M}$ be $N_{m}$ and let the total number
of I/O pins in $\mathcal{P}_{io}$ be $N_{io}$. Let $N=N_{m}+N_{io}$
be the total number of modules and I/O pins to be placed.

Let $z=(x,y)\in\mathbb{R}^{2N}$ with 
\begin{align}
x=\left(x_{1},\cdots,x_{N}\right),\\
y=\left(y_{1},\cdots,y_{N}\right),
\end{align}
by stacking the $x$-coordinates and then the $y$-coordinates of
modules from $\mathcal{M}$ and I/O pins from $\mathcal{P}_{io}$.
This stacking operation establishes an injective linear mapping from
the coordinates of modules and I/O pins to $\mathbb{R}^{2N}$. Thus,
we have two representations for the the coordinates of modules from
$\mathcal{M}$ and pins from $\mathcal{P}$. Both representations
are used in this work by the stacking convention. The representation in $\mathbb{R}^{2N}$ is used for establishing the
feasibility-seeking formulation and for introducing the feasibility-seeking
algorithm, while the 2-dimensional representation in the form $(x_{i},y_{i})$
is used for implementation.

For the \textbf{Boundary} condition I, for $1\le i\le N_{m}$, let
\begin{align}
B_{i}^{x}(z) & :=\{z\in\mathbb{R}^{2N}|\,0\leq x_{i}\leq W-w_{i}\},\\
B_{i}^{y}(z) & :=\{z\in\mathbb{R}^{2N}|\,0\leq y_{i}\leq H-h_{i}\}.\label{eq:B_i^x}
\end{align}
If both module $m_{i}$ and module $m_{j}$ fall in the floorplanning region, the following must hold, 
\begin{align}
z\in B_{i,j} & :=B_{i}^{x}\cap B_{i}^{y}\cap B_{j}^{x}\cap B_{j}^{y}.\label{box-constraints}
\end{align}
for $1\le i,j\le N_{m}$ and $i\neq j$.

The \textbf{Non-overlap} condition that the module $m_{i}$
and the module $m_{j}$ have no overlap, for $1\le i,j\le N_{m}$ and
$i\neq j$ is equivalent to one of the following four constraints
\begin{align}
z \in O_{i,j}^x & \Longleftrightarrow\mbox{\ensuremath{m_{i}} is to the left of \ensuremath{m_{j}}},\label{eq:vi:above:vj}\\
z \in O_{j,i}^x & \Longleftrightarrow\mbox{\ensuremath{m_{i}} is to the right of \ensuremath{m_{j}}},\\
z \in O_{i,j}^y & \Longleftrightarrow\mbox{\ensuremath{m_{i}} is below \ensuremath{m_{j}}},\\
z \in O_{j,i}^y & \Longleftrightarrow\mbox{\ensuremath{m_{i}} is above \ensuremath{m_{j}}},
\end{align}
where
\begin{align}
 & O_{i,j}^{x}(z):=\{z\in\mathbb{R}^{2N}|\,x_{i}+w_{i}\leq x_{j}\},\\
 & O_{i,j}^{y}(z):=\{z\in\mathbb{R}^{2N}|\,y_{i}+h_{i}\leq y_{j}\}.\label{eq:O_ij^x}
\end{align}
Then the \textbf{Non-overlap} condition is equivalent to the following
constraint 
\begin{align}
z\in O_{i,j} & :=O_{i,j}^{x}\cup O_{j,i}^{x}\cup O_{i,j}^{y}\cup O_{j,i}^{y},\label{non-overlapping-constraints}
\end{align}
for $1\le i,j\le N_{m}$ and $i\neq j$.

Let 
\begin{align}
C_{i,j} & :=O_{i,j}\cap B_{i,j}\\
 & :=\left(O_{i,j}^{x}\cap B_{i,j}\right)\cup\left(O_{j,i}^{x}\cap B_{i,j}\right)\\
 & \cup\left(O_{i,j}^{y}\cap B_{i,j}\right)\cup\left(O_{j,i}^{y}\cap B_{i,j}\right)\\
 & :=C_{i,j,L}\cup C_{i,j,R}\cup C_{i,j,B}\cup C_{i,j,A},\label{eq:union-convex-set}
\end{align}
where $L,R,B\:\mathrm{and\:}A$
stand for the relative relationship of the two modules, i.e., left, right, below, and above, respectively.
Combining the constraints in \eqref{box-constraints} and \eqref{non-overlapping-constraints},
of both the \textbf{Boundary} and the \textbf{Non-overlap} conditions,
the \textbf{floorplanning} becomes the feasibility-seeking problem:
Find a point $z$ such that
\begin{align}
z\in\bigcap_{1\leq i<j\leq N_{m}}C_{i,j}.\label{eq:feas-seek-basic-floorplanning}
\end{align}
This is a typical feasibility-seeking problem of finding a point in
the intersection of a number of sets, see, e.g., \cite{Bauschke1996}. The floorplanning
problem formulated as a feasibility-seeking problem does not handle
the \textbf{Wirelength} condition. In Section~\ref{superiorization},
below, we employ the superiorization method for handling the \textbf{Wirelength}
condition.

Furthermore, if the positions of I/O pins allow changes along the boundary
then extra constraints must be imposed. 
\begin{enumerate}
\item[IV] \textbf{I/O assignment}: Assign I/O pins to the corresponding boundaries
of the floorplanning region. 
\end{enumerate}
If the I/O pin $P_{p_{i}}$ is at the left boundary of floorplanning region,
then 
\begin{align}
D_{p_{i}}^L(z):=\left\{ z\in\mathbb{R}^{2N}\vert x_{p_{i}}^{pin}=0\;\textrm{and}\;0\leq y_{p_{i}}^{pin}\leq H\right\} .
\end{align}
Similar constraints can be constructed for I/Os at the right, top
and bottom boundaries of the floorplanning region. $\mathcal{P}_{io}\subset\mathcal{P}$
is the set, of size $N_{io}$, of I/O pins, as mentioned above. The
\textbf{FSP model of floorplanning with I/O assignment} is: 
\begin{equation}
\mathrm{Find}\text{ }z\in(\bigcap_{1\leq i<j\leq N_{m}}C_{i,j})\cap(\bigcap_{P_{p}\in\mathcal{P}_{io}}D_{p}).\label{eq:feas-seek-floorplanning-IO-assignment}
\end{equation}

\subsection{Projection Methods for Feasibility-Seeking Problems}

\label{sec:proj-for-fsp}

One method for solving feasibility-seeking problems is by using sequential
projections iteratively onto the individual sets of the family of
constraints in a predetermined order. The projection of a point $z\in\mathbb{R}^{2N}$
onto a set $C\subset\mathbb{R}^{2N}$ is the set-valued mapping $P_{C}$
\begin{equation}
P_{C}(z)=\mathop{\arg\min}\limits _{c\in C}\|z-c\|,
\end{equation}
where $\|\cdot\|$ is the Euclidean norm. When $C$ is convex, $P_{C}(z)$
is a singleton. When $C$ is non-convex, $P_{C}(z)$ may contain multiple
points, which is the situation in our case, as described below. This is the \textbf{method of alternating projection (MAP)}, see, e.g.,\cite{Escalante2011}.

Although each set $C_{i,j}$ in \eqref{eq:feas-seek-basic-floorplanning}
is not convex, the MAP method can still be applied to \eqref{eq:feas-seek-basic-floorplanning}
by selecting appropriate initial guess to guarantee the algorithmic
convergence. We do this although we are not yet able to verify that the conditions imposed in \cite{chretien1996cyclic} hold here. Nevertheless, we need to determine
which point in $P_{C_{i,j}}(z)$ to choose to enter the next iteration
in the implementation. Moreover, the order for choosing which set
$C_{i,j}$ to project onto also affect the convergence and performance
of the MAP. The algorithmic details are presented and discussed in
Section \ref{chap:algorithm}.

\section{Proposal of a Perturbed Resettable Method of Alternating Projection (Per-RMAP)}

\label{chap:algorithm}

\begin{figure}[H]
\centering \includegraphics[scale=0.3]{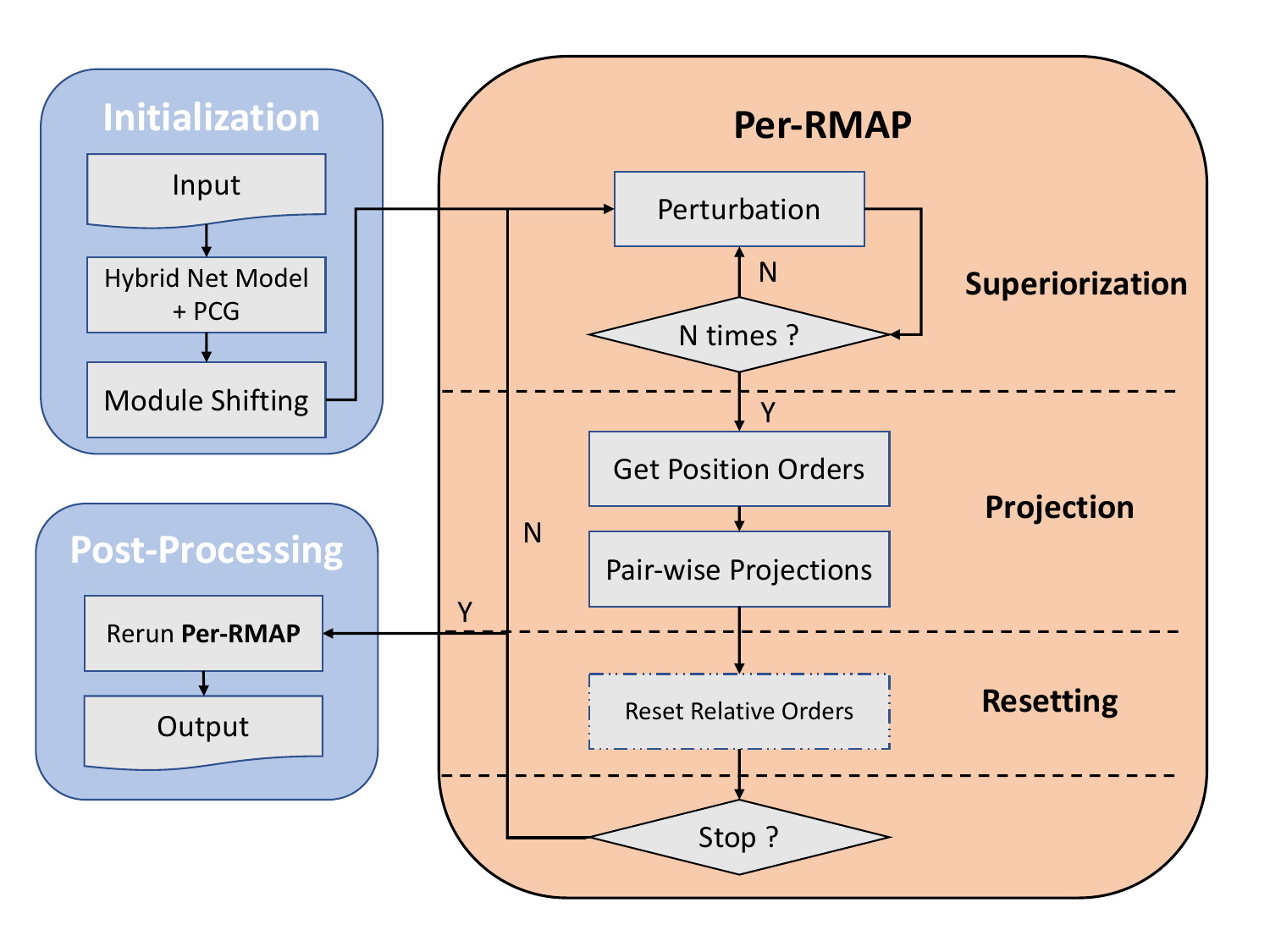}
\caption{Algorithmic Flow}
\label{fig:algorithm-framwork}
\end{figure}

In this section, the algorithmic flow is presented. As shown in Fig.~\ref{fig:algorithm-framwork},
it consists of three phases: initialization (Section~\ref{sec:algo-init}),
global floorplanning (i.e., Per-RMAP in Section~\ref{sec:algo-per-rmap}),
and post-processing (Section~\ref{sec:algo-post-proc}).

\subsection{Initialization}

\label{sec:algo-init}

The initialization can influence the final result. There are two steps
in the process for initialization of module positions. The first step
assigns modules to the positions with minimized HPWL. The second step
adjusts positions of some modules as they are easily influenced by
other modules.

In the first step, one totally ignores the cell overlaps and solves
a quadratic programming by the preconditioned conjugate gradients method (PCG) to get a total wirelength minimization, which
is similar to other force-directed placers such as SimPL~\cite{Kim2012}
and FastPlace~\cite{Viswanathan2005}. The hybrid net model \cite{Viswanathan2005}
is used for net decomposition. It is a combination of a classical
clique model and a star model, which not only gets a trade-off between
speed and accuracy but also better captures the relative orders after
net decomposition. However, the modules tend to cluster together if
there is little connection between modules and I/O pins that lie on
the boundary of the floorplanning region.

In the second step, modules shifting is used to improve the initialization.
That is, detect key modules that may impact the final result and initialize them to the boundary of the floorplanning region. 

\subsection{Global Floorplanning: Per-RMAP}

\label{sec:algo-global-fplan} In global floorplanning we design a
Per-RMAP algorithm based on projections, where a resetting strategy
and the superiorization methodology are utilized. 

\subsubsection{Resettable Method of Alternating Projection (RMAP)}

While MAP (introduced in Section~\ref{sec:proj-for-fsp}) solves
an FSP with convex constraints sets, resetting strategy is designed
to tackle the situation when the constraints sets are unions of convex
sets. To avoid getting stuck at an infeasible solution or oscillating
among infeasible solutions, our idea is to generalize MAP into a ``strategy-enabled MAP''. The insight is that MAP always chooses the subset in the union
closest to the current iteration, and some choices prevent a feasible
solution. We enable a customizable choice strategy via the ``preference
ratio'' among the subsets in the union. The details are shown in
Algorithm~\ref{alg:Framwork}.
\begin{algorithm}[hbt]
\caption{Generalized Projection with Preferences (it becomes RMAP when $\mathit{preference\_ratio}$ adopts the resetting strategy)}
\label{alg:Framwork} 
\begin{algorithmic}[1]
\Require
  \Statex The initial positions: $z=(x,y)\in\mathbb{R}^{2N}$
  \Statex  The processing order of constraints: $order$
\Ensure
  \Statex The updated positions: $z=(x,y)$
  \Statex \vspace{-0.5em}
\For{$C_{i,j}$ in $order$}
  \State $(\eta_{L},\eta_{R},\eta_{B},\eta_{A}) = preference\_ratio(i,j)$
  \State $w_{t} = {\exp({\eta_{t}}/{\epsilon})}/ {\sum_{k}{\exp({\eta_{k}}/{\epsilon})}}$, for $t\in\{L,R,B,A\}$
  \State $z=\sum_{t\in\left\{L,R,B,A\right\}}w_{t} \cdot P_{C_{i,j,t}}(z)$
\EndFor 
\State \Return $z$
\end{algorithmic} 
\end{algorithm}
The iterative process scans modules in a certain order and applies pair-wise projections which are the weighted average
of 4 projections onto convex sets. A key issue is the correct choice
of the convex set to project onto, i.e., the distribution of the preference
ratios. The preference ratio gives the preference of the choice of
the four convex sets and is amplified by an exponential function.
The MAP uses the closest point strategy for the preference ratio,
which is calculated by $\eta_{t}=-\|z-P_{i,j,t}(z)\|$. A resetting
strategy based on the closest point strategy is designed, which considers
previous behaviors of projection for the preference ratio calculation,
as shown in \eqref{eq:resetting-strategy}, 
\begin{equation}
\eta_{k}=\left\{ \begin{array}{lll}
-\infty, & c_{k}>T & \text{(and reset $c_{k}=0$)},\\
-\|z-P_{k}(z)\|, & otherwise & \text{(and set $c_{k}=c_{k}+1$)},
\end{array}\right.\label{eq:resetting-strategy}
\end{equation}
where $T$ is a predefined positive integer. The $c_{k}$ is the count
of each convex set $C_{i,j,k},\textrm{ for }k\in\{L,R,B,A\}$ that has been projected onto since the last reset.

That is, when a pair of modules repeatedly (for more than $T$ times) projects in a certain direction but fails to remove
overlap, a ``reset'' action is activated, and this direction is given a lowest preference ratio in this iteration to escape from the oscillation or the stuck situation.

The processing order can be the position order, such that modules with smaller $x$ coordinates and  $y$ coordinates are processed earlier. It could
preserve the relative order of most modules.

\subsubsection{The Superiorization Methodology (SM)}

\label{superiorization} The SM lies between feasibility-seeking and
constrained optimization.  While seeking compatibility with constraints,
SM reduces the value of an objective function but not necessarily
to a minimum. It takes proactive measures to perturb its iterates
in a manner that guides them towards a feasible point and reduces
the value of the objective function. In our work, we report an improvement
in the total wirelength when applying novel modifications of superiorization.
Inspired by a modern version of superiorization~\cite{yair2020},
our superiorization for HPWL improvement is shown in Algorithm~\ref{alg:superiorization}.
\begin{algorithm}[hbt]
\caption{Superiorization Methodology (SM)}
\label{alg:superiorization} \begin{algorithmic}[1] \Require \Statex
The intermediate positions: $z=(x,y)\in\mathbb{R}^{2N}$ 
 \Statex Number of perturbations in one iteration: $\mbox{\textit{Num}}$ 
 \Statex Current iteration number: $k$ 
 \Statex Perturbation decay index: $\ell_{k-1}$ 
 \Statex Minimum perturbation length: $\lambda_{min}$ 
 \Statex Initial perturbation length: $\lambda_{init}$
 \Statex Perturbation decay factor: $\Lambda\in(0,1)$ 
 \Ensure \Statex The updated
positions: $z$ \Statex The new perturbation decay index: $\ell_{k}$
\Statex \vspace{-0.5em}
\For{$n=1:\mbox{\textit{Num}}$}\label{eq:pert-loop-1} 
    \If{$k<l_{k-1}$}\label{eq:pert-decay-index}
    \State $\ell_{k} = \text{a random integer in $\left[k,\ell_{k-1}\right]$}$
    \Else \State $\ell_{k}=k$ 
    \EndIf
 \State $v^{k,n}=\nabla \mbox{\textrm{HPWL}}(x,y)$ 
 \For{$cnt=1:10$} 
 \State $\lambda_{pert}^{k,n}=max(\lambda_{min},\lambda_{init}\cdot\Lambda^{l_{k}})$\label{eq:lambda-pert} 
\State $(x',y')=z'=z-\lambda_{pert}^{k,n}\cdot{v^{k,n}}/{\|v^{k,n}\|}$\label{eq:pos-pert}
\If{$\mbox{\textrm{HPWL}}(x',y')<\mbox{\textrm{HPWL}}(x,y)$} 
\State $z=z'$ 
\State \textbf{break} // and continue at line \ref{eq:pert-decay-index}
\EndIf
\State $\ell_{k}=\ell_{k}+1$ 
\EndFor \EndFor \State \Return
$(z,\ell_{k})$ 
 \label{eq:pert-loop-2} 
\end{algorithmic} 
\end{algorithm}
At iteration $k$, perturbations are applied $\mbox{\textit{Num}}$ times to move
the modules and I/O pins along the negative gradient of the HPWL function.
The step-size parameter of perturbation $\lambda_{pert}^{k,n}$ (line
\ref{eq:lambda-pert}) decreases with iterations. We design the step-size
parameter $\lambda_{pert}^{k,n}$ by considering the following properties: 
\begin{enumerate}
\item Step sizes of the perturbations $\lambda_{pert}^{k,n}$ must be summable,
i.e., $\sum_{k=0}^{\infty}\sum_{n=0}^{\mbox{\textit{Num}}-1}\lambda_{pert}^{k,n}<+\infty$.
In our algorithm, the step sizes are generated via a subsequence of
 $\{\Lambda^{\ell}\}_{\ell=0}^{\infty}$ with $\ell$ powers of the user-chosen kernel
$0<\Lambda<1$. 
\item A lower bound $\lambda_{min}$ is given to ensure the performance
of the perturbation. In our setting, $\lambda_{min}=0.1$. 
\item Controlling the decrease of the step sizes of objective function reduction perturbations in the superiorization. If the step sizes $\lambda_{pert}^{k,n}$ decrease too fast (i.e., $\ell_{k-1} > k$ in $\Lambda^{\ell_{k-1}}$ from the previous SM iteration), then too little leverage is allocated to the total wirelength reduction. So the perturbation decay index $\ell_k$ is a
number between the current iteration index $k$ and the value of $\ell_{k-1}$ from the last iteration sweep, as shown in line \ref{eq:pert-decay-index}. This modification was suggested and investigated in \cite{Langthaler2014} \cite{Prommegger2014}.
\end{enumerate}
If the perturbation (line \ref{eq:lambda-pert}-\ref{eq:pos-pert}) shortens the total wirelength, then the algorithm accepts the perturbation;
otherwise, we repeatedly decrease the step size for at most 10 times
until a shorter wirelength is reached.

\subsubsection{Perturbed Resettable Method of Alternating Projection (Per-RMAP)}\label{sec:algo-per-rmap} 

Combining RMAP and SM, Per-RMAP is designed, as shown in Algorithm~\ref{alg:PerRMAP}.
 
\begin{algorithm}[hbt]
\caption{ Perturbed Resettable Method of Alternating Projection (Per-RMAP)}
\label{alg:PerRMAP} \begin{algorithmic}[1] \Require 
\Statex The initial positions: $z=(x,y)\in\mathbb{R}^{2N}$ 
\Statex Number of perturbations in one iteration: $\mbox{\textit{Num}}$ 
\Statex Minimum perturbation length: $\lambda_{min}$ 
\Statex Initial
perturbation length: $\lambda_{init}$ 
\Statex Perturbation decay factor: $\Lambda\in(0,1)$ 
\Statex Initial projection length: $\gamma_{init} \in (0,1)$
\Statex Projection progress factor: $\Gamma>1$
\Ensure \Statex The updated positions $z=(x,y)$ \Statex \vspace{-0.5em}
\State $\ell_{-1}=0$
\For{$k=0:\infty$}
\State $(z,\ell_{k})=\mbox{\textrm{SM}}(z,\mbox{\textit{Num}},k,\ell_{k-1},\lambda_{min},\lambda_{init},\Lambda)$
\State $order=$ {[}generated by position order{]}
\State $\gamma_{proj}=min(1,\gamma_{init}\cdot\Gamma^{k})$\label{eq:lambda-proj}
 \State $z=z+\gamma_{proj}\cdot(\mbox{\textrm{RMAP}}(z,order)-z)$ 
 \State $isStop = \mbox{\textrm{RelativeOverlappingAreaCheck}}(z)$
 \If{$isStop == true$}
 \State \textbf{break}
 \EndIf
\EndFor
\State \Return $z$
\end{algorithmic} 
\end{algorithm}

In Per-RMAP, superiorization with $\mbox{\textit{Num}}$ perturbations is firstly
applied, projections with the resetting strategy in position order are followed.
The relaxation parameter of the projections (line \ref{eq:lambda-proj})
increases with iterations and has an upper bound. As a result, the
projection steps become a dominant part as iterations proceed and
does not move the cells too far in any one iteration. 

\subsection{Post Processing}\label{sec:algo-post-proc}
After global floorplanning, we obtain a result with relative overlapping area less than $0.1\%$. The
superiorization, which may drag modules closer, has less and less
impact on the position changing at iterations. As a consequence, there
may still exist some gaps between modules. Hence RMAP is rerun to
improve the result. At this phase, perturbation decay index is reset
to $k\times\epsilon$, where $k$ is the total iteration number and
$\epsilon\in(0,1)$.

\section{Experiment}
\subsection{Benchmark}
The Microelectronics Center of North Carolina (MCNC)
benchmark \cite{Kozminski1991} is a commonly used benchmark, which consist of five instances. The details are shown in Table~\ref{tab:MCNC-1}, including the number of modules, I/O pins, pins, nets as well as the size of the floorplanning region (die size).
\begin{table}[H]
\centering{}\caption{Characteristics of the MCNC block instances.}
\label{tab:MCNC-1} %
\begin{tabular}{cccccc}
\toprule
\multirow{2}{*}{Instance}  & \multicolumn{4}{c}{Number of} &\multirow{2}{*}{Die Size}\\
\cmidrule(lr){2-5}
 & Modules  & I/O pins & Pins  & Nets &\\
\midrule
\textbf{apte}  & 9  & 73  & 214  & 97 & 10,500$\times$10,500\\
\textbf{xerox}  & 10  & 2  & 696  & 203 & 5831$\times$6412\\
\textbf{hp}  & 11  & 45  & 264  & 83 & 4928$\times$4200\\
\textbf{ami33}  & 33  & 42  & 480  & 123 & 2058$\times$1463\\
\textbf{ami49}  & 49  & 22  & 931  & 408 & 7672$\times$7840\\
\bottomrule
\end{tabular}
\end{table}

\subsection{Performance of The Resetting Strategy}

Table \ref{tab:compare-MAP-and-RMAP} shows the performance of the reset
strategy in our experiment on the MCNC benchmarks. The procedure stops
when the \textbf{relative overlapping area (Relative O. A.)} oscillates
among some values, stays constant at a value or is less than $0.1\%$.
Results on MCNC benchmarks show that without the resetting strategy, the
MAP may get stuck at infeasible points, where the overlaps are not
totally removed. However, with the resetting strategy, this phenomenon is removed.

\begin{table}[ht]
    \caption{MAP versus RMAP on the MCNC benchmarks. Results show that RMAP can relieve oscillations.}
    \label{tab:compare-MAP-and-RMAP}
    \centering
    \begin{tabular}{ccccc}
        \toprule
        Instance & Reset? & Runtime & Iterations & Relative O. A.\\
        \midrule
        \multirow{2}{*}{apte} & N & 0.41 & 158 & 11.5\%\\        
         & Y & 0.35 & 41 & $< 0.1\%$\\
        \midrule
        \multirow{2}{*}{xerox} & N & 0.27  & 131 & 10.5\%\\
         & Y & 0.13 & 36 & $< 0.1\%$\\
        \midrule
        \multirow{2}{*}{hp} & N & 0.47 & 286 & 0.6\%\\
         & Y & 0.06 & 22 & $< 0.1\%$\\
        \midrule
        \multirow{2}{*}{ami33} & N & 1.11 & 279 & 2.4\%\\
         & Y & 0.71 & 80 & $< 0.1\%$\\
        \midrule
        \multirow{2}{*}{ami49} & N & 8.13 & 931 & 3.4\%\\
         & Y & 3.25 & 215 & $< 0.1\%$\\
        \bottomrule
    \end{tabular}
\end{table}
\subsection{Floorplanning Results}

We compare our result with some state-of-the-art results \cite{funke2016exact} obtained by a branch-and-bound (B\&B) method. This method achieves the optimal floorplans on the first three instances (apte, xerox, and hp). However, it only obtains sub-optimal floorplans on the last two larger instances (ami33 and ami49) due to time limit. At the initialization step, key modules of smaller size compared with other modules in the same instance are initialized to the boundary of the floorplanning region\footnote{For the instance xerox, it is the module named $\mbox{\textit{BLKP}}$. For the instance hp, it is modules named $cc\_22$ and $new1$.}. For original floorplanning without I/O assignment,
results are shown in Table \ref{tab:floorplan-result-wo-io-1}. The
hyper-parameters are $\lambda_{min}=0.1, \Lambda=0.99,\lambda_{init}=321,\gamma_{init}=0.7804,\Gamma=1.1,\epsilon=0.35$. The HPWL of our method is only $3\%$ longer than the optimal results in B\&B~\cite{funke2016exact} within $1.95\times$ execution time\footnote{Our prototype is in MATLAB and is expected to be fast if using C++.}.

\begin{table}[H]
\caption{Experimental Results for Floorplanning without I/O assignment.}
\label{tab:floorplan-result-wo-io-1} \centering \setlength{\tabcolsep}{3pt}
\begin{threeparttable} %
\begin{tabular}{lllllll}
\toprule
\multirow{2}{*}{Instance } & \multicolumn{2}{l}{\multirow{2}{*}{B\&B~\cite{funke2016exact}}} & \multicolumn{2}{l}{Per-RMAP} & \multicolumn{2}{l}{Ratio}\\
 & & & \multicolumn{2}{l}{w.o. I/O Assignment} & \multicolumn{2}{l}{Per-RMAP/B\&B}\\
\cmidrule(r){2-3}\cmidrule(r){4-5}\cmidrule(r){6-7}
 & HPWL  & Time (sec)  & HPWL  & Time (sec)  & HPWL  & Time \\
\midrule
\textbf{apte}  & 513061  & 13  & 528618  & 8.84  & 1.03  & 0.68\\
\textbf{xerox}  & 370993  & 48  & 382596  & 172.09  & 1.03  & 3.59 \\
\textbf{hp}  & 153328  & 102  & 159979  & 5.06  & 1.04  & 0.05\\
\textbf{ami33}  & 58627  & 13  & 61444  & 23.83  & 1.05  & 1.83\\
\textbf{ami49}  & 640509  & 73  & 637098  & 261.76  & 1.00  & 3.59\\
\textbf{Average}  & -  & -  & -  & -  & 1.03  & 1.95 \\
\bottomrule
\end{tabular}\end{threeparttable} 
\end{table}
For original floorplanning with I/O assignment, results are shown
in Table \ref{tab:floorplan-result-with-io}. The hyper-parameters
are $\lambda_{min}=0.1, \Lambda=0.99$, $\lambda_{init}=488$, $\gamma_{init}=0.7761$, $\Gamma=1.0001$, and $\epsilon=0.35$.
Compared with B\&B, our method achieves 8$\%$ improvement on total wirelength and takes $4.99\times$ time over B\&B \cite{funke2016exact}.

\begin{table}[H]
\caption{Experimental Results for Floorplanning with I/O Assignment.}
\label{tab:floorplan-result-with-io} \centering \setlength{\tabcolsep}{3pt}
\begin{threeparttable} %
\begin{tabular}{lllllll}
\toprule
\multirow{2}{*}{Instance} & \multicolumn{2}{l}{\multirow{2}{*}{B\&B~\cite{funke2016exact}}} & \multicolumn{2}{l}{Per-RMAP} & \multicolumn{2}{l}{Ratio}\\
 & & & \multicolumn{2}{l}{with I/O Assignment} & \multicolumn{2}{l}{Per-RMAP/B\&B}\\
\cmidrule(r){2-3}\cmidrule(r){4-5}\cmidrule(r){6-7}
 & HPWL  & Time (sec)  & HPWL  & Time (sec)  & HPWL  & Time \\
\hline 
\textbf{apte}  & 513061  & 13  & 362587  & 122  & 0.71  & 9.39\\
\textbf{xerox}  & 370993  & 48  & 382587  & 188  & 1.03  & 3.92 \\
\textbf{hp}  & 153328  & 102  & 149545  & 139  & 0.98  & 1.36\\
\textbf{ami33}  & 58627  & 13  & 54489  & 79  & 0.93  & 6.08\\
\textbf{ami49}  & 640509  & 73  & 618898  & 306  & 0.97  & 4.19\\
\textbf{Average}  & -  & -  & -  & -  & 0.92  & 4.99 \\
\bottomrule
\end{tabular}\end{threeparttable} 
\end{table}
\label{chap:experiment}

\section{Conclusion and Future Work}

\label{chap:conclusion} We model the fixed-outline floorplanning
as a feasibility-seeking problem (FSP). However, the conventional
method of alternating projection (MAP) for FSP cannot obtain legal
floorplans. This is because the constraints sets of the floorplanning problem are not convex. We analyze
the union convex property of the constraints sets and propose the resettable method of alternating projection (RMAP) to improve its convergence to a feasible solution. Furthermore, a superiorized version, Per-RMAP, is designed to decrease the total wirelength. The experiments show that our method achieves nearly optimal results for some floorplanning problems and 8$\%$ improvement compared with branch-and-bound (B\&B) method on half-perimeter wirelength (HPWL) after considering I/O assignment. Our future work is to investigate the scalability of our method by adding more experiments on larger instances and also investigate the capability of handling complex practical constraints, like the ones for 2.5D floorplanning.

\section*{Acknowledgment}
We thank our colleagues Tie Zhou and Jiansheng Yang at the Peking
University for their helpful comments during our joint Zoom meetings
on this project. 

This research is supported by National Natural Science Foundation of China (NSFC) under Grant No. 11961141007 and by the Israeli Science Foundation (ISF) under Grant No. 2874/19 within the NSFC-ISF joint research program. The work of Y.C. is also supported by the U.S. National Institutes of Health Grant No. R01CA266467. Additionally, the work of G.L. is supported by National Key R\&D Program of China under Grant No. 2022YFB4500500.

\balance

\bibliographystyle{IEEEtran}
\bibliography{ISEDA-FSP}

\end{CJK*}
\end{document}